\documentclass[12pt,a4paper]{article}
\usepackage{amsmath,amsthm,amsfonts,latexsym,amscd,amssymb, makeidx,graphics, caption}
\makeindex
\usepackage{float}
\usepackage{enumitem}
\usepackage{tabularx}
\usepackage{breqn}
\newtheorem{theorem}{Theorem}[section]
\newtheorem{lemma}[theorem]{Lemma}
\newtheorem{corollary}[theorem]{Corollary}
\newtheorem{proposition}[theorem]{Proposition}

\theoremstyle{definition}

\theoremstyle{remark}
\newtheorem*{remark}{Remark}

\usepackage{mathtools}

\def\addsymbol #1: #2#3{$#1$ \> \parbox{5in}{#2 \dotfill \pageref{#3}}\\}

\begin{document}

	\begin{center}
		\Large{The Structure of Orthomorphism Graph of $\mathbb{Z}_2 \times \mathbb{Z}_4$} 
	\end{center}
	
	\begin{center}
		Rohitesh Pradhan and Vivek Kumar Jain\\
		Central University of South Bihar, Gaya, India\\
		Email:rohiteshpradhan@gmail.com, jaijinenedra@gmail.com
	\end{center}
	\large{\textbf{Abstract:}}
	In this paper, we gave a theoretical proof of the fact that Orthomorphism graph of group $\mathbb{Z}_2 \times \mathbb{Z}_4$ has maximal clique 2, by determining the structure of the graph.
	\vspace{.5cm}
	
	\noindent \textbf{\textit{Keywords:}} Orthomorphism Graph
	\smallskip
	
	\noindent \textbf{\textit{2020 Mathematical Subject Classification:}} 05B15
	\section{Introduction}

	Let $G$ be a group. A bijection $\theta$: $G \rightarrow G$ for which map $\phi$: $x \mapsto x^{-1}\theta(x)$ is also a bijection of $G$ is called \textit{Orthomorphism of G}. Two orthomorphisms $\theta_{1}$ and $\theta_{2}$ of $G$ are called \textit{orthogonal}, written  $\theta_{1} \perp \theta_{2}$, if $x\mapsto\theta_{1}(x)^{-1}\theta_{2}(x)$ is  a bijection of $G$. An orthomorphism of a group which fixes identity element of the group is called \textit{normalised} orthomorphism. Now onwards by an orthomorphism we mean normalised orthomorphism. We denote the set of orthomorphism of a group $G$ by \textit{Orth($G$)}. A graph in which vertices are orthomorphisms of $G$ and adjacency being synonymous with orthogonality is called \textit{orthomorphism graph} of $G$, which is also denoted by Orth($G$). The order of the largest complete subgraph of a graph is called \textit{clique number} of the graph. Clique number of orthomorphism graph of a group $G$ is denoted by $\omega(G)$.
	In this paper, we prove that $\omega(\mathbb{Z}_2 \times \mathbb{Z}_4) = 2$.
	
 In 1961, Johnson, Dulmage, and Mendelsohn showed that $\vert Orth(\mathbb{Z}_2 \times \mathbb{Z}_4) \vert$ = 48 and $\omega(\mathbb{Z}_2 \times \mathbb{Z}_4) = 2$, via a computer search. In 1964, through exhaustive hand computation, Chang, Hsiang, and Tai \cite{chang} also found that $\omega(\mathbb{Z}_2 \times \mathbb{Z}_4)$ = 2 and in 1986, via a computer search, Jungnickel and Grams \cite{grams} also confirm above fact. In 1992, Evans and Perkel found that $Orth(\mathbb{Z}_2 \times \mathbb{Z}_4)$ consists of 12 disjoint 4-cycles using Cayley (a forerunner of the computer algebra system Magma)\cite{evans} and asked for theoretical proof of this fact \cite[Problem 19]{evans}. In 2021, Evans gave a theoretical proof of this in \cite{evans3}. In this paper we gave another proof of this fact.
	
	An \textit{automorphisms} $A$ of Orth($G$) is a bijection on Orth($G$) such that $A(\theta_{1}) \perp A(\theta_{2})$ if and only if $\theta_{1} \perp \theta_{2}$ where $\theta_{1}, \theta_{2} \in$ Orth($G$).
	
	For $f \in$ Aut($G$), the group of automorphism of $G$, the map $H_{f}\colon$ Orth($G$) $\rightarrow$ Orth($G$) defined as $H_{f}(\theta) = f\theta f^{-1}$ is known as \textit{homology} of Orth($G$). Homology is an example of automorphism of Orth($G$) \cite[Theorem 8.6, p.206]{evans2}. Any unexplained notation used in the paper is from \cite{evans2}.  
	\section{Basic Results}
	
	Suppose $G=\{g_{1}=e, g_{2}, \ldots,g_{n}\}$ is a finite group. $e$ is identity element of $G$. Suppose $\tau$ denotes the regular left permutation representation of $G$. Let us identify G with $\tau(G)$. So $G \leqslant Sym(G)$. Further, let us identify $Sym(G)$ with $Sym\{{1,2,\ldots,n}\}$ (denoted as $S_{n}$) by identifying $g_{i}$ with $i$. Clearly, if $\theta \in$ Orth($G$), then $\theta \in$ Sym\{2,3,$\ldots$,n\}(denoted as $S_{n-1}$). For $\theta \in$ Orth($G$), the map $\phi_{\theta}\colon$ $G \rightarrow G$ defined as $\phi_{\theta}(x) = x^{-1}\theta(x)$ is called the \textit{complete mapping} associated with $\theta$. Clearly $\phi_{\theta} \in S_{n}$. Also note that $\theta(x) = x$ if and only if $x$ is identity element of $G$.
	\begin{remark}
		For $\theta \in Orth(\mathbb{Z}_2 \times \mathbb{Z}_4)$, take $a, y \in \mathbb{Z}_{2} \times \mathbb{Z}_{4}$ such that $o(a)=4$ and $o(y) = o(\theta(y)) = 2$. Then $\{x \in \mathbb{Z}_2 \times \mathbb{Z}_4 \mid o(x) = 4\} = \{a, ay, a\theta(y), ay\theta(y)\}$ and $ \{x \in \mathbb{Z}_2 \times \mathbb{Z}_4 \mid o(x) = 2\} = \{y, \theta(y), y\theta(y)\}.$
	\end{remark}
	\begin{lemma}\label{le1}
	Let $G=\mathbb{Z}_{2} \times \mathbb{Z}_{4}$ be a group and $\theta$ is a map from $G$ to $G$. Define $A_{ij}$ $i,j \in \{2,4\}$ as follows \begin{align*}
	A_{44} =& \{x \in G \mid o(x) = 4, o(\theta(x)) = 4 \},\\
	A_{42} =& \{x \in G \mid o(x) = 4, o(\theta(x)) = 2 \},\\
	A_{24} =& \{x \in G \mid o(x) = 2, o(\theta(x)) = 4 \},\\
	A_{22} =& \{x \in G \mid o(x) = 2, o(\theta(x)) = 2 \}.
\end{align*}
\begin{enumerate}[wide, labelindent=*]
	\item[(1)]  $\theta$ is a bijection fixing identity element of $G$ if and only if
\begin{enumerate}
	\item [(a)] $\{\theta(x)\mid x \in A_{44}\} \sqcup \{\theta(x)\mid x \in A_{24}\}$ = $\{x \in G \mid o(x)= 4\}$, and
	\item[(b)] $\{\theta(x)\mid x \in A_{42}\} \sqcup \{\theta(x)\mid x \in A_{22}\}$ = $\{x \in G \mid o(x)= 2\}$.
\end{enumerate}
\item[(2)] $\theta$ is an orthomorphism if and only if $\theta$ is   bijection fixing identity and
\begin{enumerate}
	\item [(a)]$\{x^{-1}\theta(x)\mid x \in A_{42}\} \sqcup \{x^{-1}\theta(x)\mid x \in A_{24}\}$ = $\{x \in G \mid o(x)= 4\}$, and
	\item [(b)]$\{x^{-1}\theta(x)\mid x \in A_{44}\} \sqcup \{x^{-1}\theta(x)\mid x \in A_{22}\}$ = $\{x \in G \mid o(x)= 2\}$.
\end{enumerate}
\end{enumerate}
	\end{lemma}
	\begin{proof}
		Follows from the definition of bijection and orthomorphism.
	\end{proof}
\begin{corollary}\label{co1}
If $\theta \in Orth(\mathbb{Z}_2 \times \mathbb{Z}_4)$, then $\vert A_{44} \vert = \vert A_{42} \vert = \vert A_{24} \vert = 2$ and $\vert A_{22} \vert = 1$. 
\end{corollary}
	\begin{corollary}\label{co2}
		For $\theta \in Orth(\mathbb{Z}_2 \times \mathbb{Z}_4)$
		\begin{enumerate}
			\item [(i)] $\vert A_{44} \cap \theta(A_{44}) \vert$ = $\vert A_{44} \cap \theta(A_{24}) \vert$ = $\vert A_{42} \cap \theta(A_{44}) \vert$ = $\vert A_{42} \cap \theta(A_{24}) \vert$ = $1$.
			\item [(ii)] If $A_{22} = \{x\}$, then $\theta(A_{42}) = \{x, x\theta(x)\}$ and $A_{24} = \{\theta(x), $ $x\theta(x)\}.$
			\item [(iii)] $\{x^{-1}\theta(x) \mid x \in A_{44}\} = \{x, \theta(x)\}$.
		\end{enumerate}
	\end{corollary}
\begin{proof}
	(i) From Corollary \ref{co1}, exactly two element of order $4$ will map to order $4$ elements. Since $a^{-1}(ay) = (ay)^{-1}a = y$, where $y$ is an element of order $2$ and $a \in A_{44}$, so $A_{44} \neq \theta(A_{44}).$ Suppose $A_{44} \cap \theta(A_{44})$ = $\phi$.
	Then if $\{a, ay\} \in A_{44}$ then $\{az, azy\} \in \theta(A_{44})$, where $z \neq y$ and $z,y$ are elements of order $2$. Then $a^{-1}(az) = (ay)^{-1}(azy) = z$ or $a^{-1}(azy) = (ay)^{-1}(az) = zy$ which contradicts the bijectivity of $\phi_{\theta}$. Hence $\vert A_{44} \cap \theta(A_{44})\vert = 1$. As $\theta(A_{44})$ is in partition with $\theta(A_{24})$ so $\vert A_{44} \cap \theta(A_{24})\vert = 1$, also $A_{44}$ and $A_{42}$ are in partition, therefore $\vert A_{42} \cap \theta(A_{44})\vert = \vert A_{42} \cap \theta(A_{24})\vert = 1$.\\
	(ii) Suppose $A_{22} = \{x\}$. Then $\theta(x) \neq x$. So $\{x, \theta(x), x\theta(x)\}$ is the set of all elements of order $2$ in $\mathbb{Z}_2 \times \mathbb{Z}_4$. $\theta(A_{42}) \sqcup \theta(A_{22})$ = $\{x,\theta(x), x\theta(x)\}$.
	Clearly, $\theta(A_{42}) = \{x, x\theta(x)\}$.
	As $x \notin A_{24}$, so $A_{24} = \{\theta(x), x\theta(x)\}$.\\
	(iii) If $\theta$ is an orthomorphism then, $\{x^{-1}\theta(x)\mid x \in A_{44}\} \sqcup \{x^{-1}\theta(x)\mid x \in A_{22}\}$ = $\{x, \theta(x), x\theta(x)\}$. Therefore $\{x^{-1}\theta(x)\mid x \in A_{44}\} = \{x , \theta(x)\}$.
\end{proof}
\begin{proposition}\label{pr1}
	If $\theta \in Orth(\mathbb{Z}_2 \times \mathbb{Z}_4)$ and $x \in A_{22}$, $a \in A_{44} \setminus \theta(A_{44})$, then $\theta(a) = ax$. Moreover $A_{44} = \{a, ax\}, \theta(A_{44}) = \{ax, ax\theta(x)\}$ and $A_{42} = \{ax\theta(x), a\theta(x)\}$.
\end{proposition}
\begin{proof}
	Suppose $x \in A_{22}$ and $a \in A_{44}\setminus\theta(A_{44})$. By Corollary \ref{co2} (iii), $a^{-1}\theta(a) \in \{x, \theta(x)\}$.  Assume $\theta(a) = a\theta(x)$. By Corollary \ref{co2} (i), $A_{44} =\{a, a\theta(x)\}$. Then by Corollary \ref{co2} (iii), $\theta(a\theta(x)) = ax\theta(x)$ and $A_{42} = \{ax, ax\theta(x)\}$. By Corollary \ref{co2} (ii), $\theta(A_{42}) = \{x, x\theta(x)\}$.\\
	\textbf{Case(a):} If $\theta(ax) = x$ and $\theta(ax\theta(x)) = x\theta(x)$, then $\phi_{\theta}(ax) = a^{-1}$ and $\phi_{\theta}(ax\theta(x)) = a^{-1}$, which contradicts the bijectivity of $\phi_{\theta}$.\\
	\textbf{Case(b):} If $\theta(ax) = x\theta(x)$ and $\theta(ax\theta(x)) = x$, then $\phi_{\theta}(ax) = a^{-1}\theta(x)$ and $\phi_{\theta}(ax\theta(x)) = a^{-1}\theta(x)$, which again contradicts the bijectivity of $\phi_{\theta}$.
	\\ Thus, $\theta(a) = ax$ and  $\theta(ax) = ax\theta(x)$. Hence $A_{44} = \{a, ax\}$, $\theta(A_{44}) = \{ax, ax\theta(x)\}$ and $A_{42} = \{ax\theta(x), a\theta(x)\}$.
\end{proof}
This can be shown by the Figure 1
\begin{figure}[H]
	\begin{center}
		\captionsetup{labelformat=empty}
	\includegraphics[scale = .43]{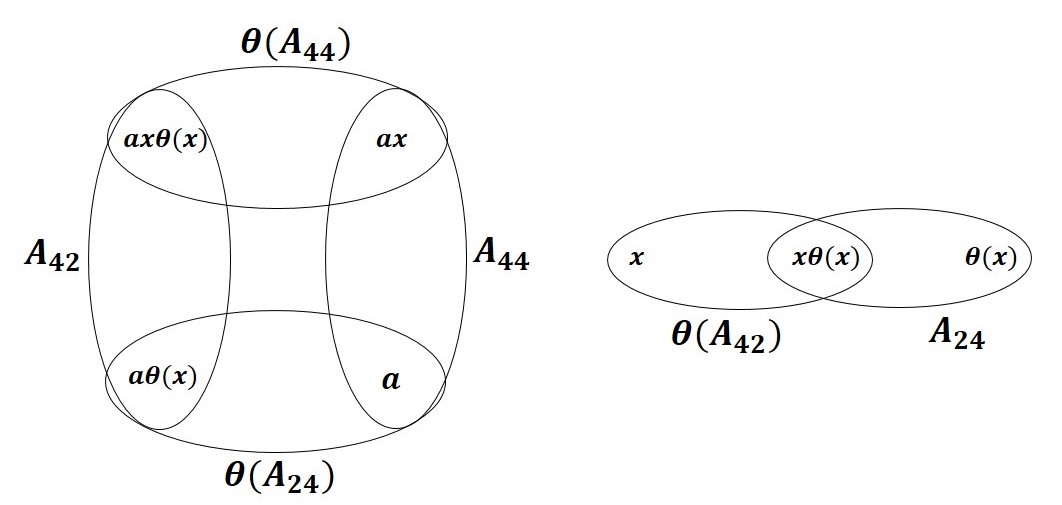} \caption{Figure 1}
\end{center}

\end{figure}	
	\section{The Structure of Orth($\mathbb{Z}_2 \times \mathbb{Z}_4$)}
\begin{theorem}\label{th1}
	Suppose $\theta \in Orth(\mathbb{Z}_2 \times \mathbb{Z}_4)$, $A_{22} = \{x\}$, $a \in A_{44} \setminus \theta(A_{44})$. Then $\theta$ has one of the following form
	\begin{enumerate}
		\item [(i)] $(a, ax, ax\theta(x), x\theta(x),a\theta(x), x, \theta(x))$, where $x\theta(x) = a^{2}$.
		\item[(ii)] $(a, ax, ax\theta(x), x\theta(x))(\theta(x), a\theta(x),x)$, where $x\theta(x) \neq a^{2}$.
		\item [(iii)] $(a,ax,ax\theta(x), x, \theta(x))(a\theta(x),x\theta(x))$, where $x\theta(x) \neq a^{2}$.
		\item[(iv)] $(a,ax, ax\theta(x),x,\theta(x),a\theta(x),x\theta(x))$, where $x\theta(x) = a^{2}$.
	\end{enumerate} 
and $\vert Orth(G) \vert = 48$.
\end{theorem}
\begin{proof}
	Suppose $A_{22} = \{x\}$ and $a \in A_{44} \setminus \theta(A_{44})$. Then by Proposition \ref{pr1}, $A_{44} =\{a, ax\}$, $\theta(A_{44}) = \{ax, ax\theta(x)\}$ and $A_{42} = \{ax\theta(x), a\theta(x)\}$.
	
	\textbf{{Case(i):}} Assume $\theta(ax\theta(x)) = x\theta(x)$.
	Clearly, $\theta(a\theta(x)) = x$. Then $\phi_{\theta}(ax\theta(x)) = a^{-1}$ and $\phi_{\theta}(a\theta(x)) = a^{-1}x\theta(x)$. 
	
	\textbf{Subcase(a):} Assume $\theta(\theta(x)) = a$. Then $\theta(x\theta(x)) = a\theta(x)$, $\phi_{\theta}(\theta(x)) = a\theta(x)$ and $\phi_{\theta}(x\theta(x)) = ax$. Bijectivity of $\phi_{\theta}$ implies, $a^{-1}x\theta(x) = a$ or $x\theta(x) = a^{2}$. Thus, if $x\theta(x) = a^{2}$, then $\theta$ is an orthomorphism, given by $(i)$. We have $4$ choices for $a$ and 2 choices for $x$ as $x \neq a^{2}$. Hence, we have 8 orthomorphisms of this form in $Orth(G)$.
	
	\textbf{Subcase(b):} Assume $\theta(\theta(x)) = a\theta(x)$. Then $\theta(x\theta(x)) = a$, $\phi_{\theta}(\theta(x)) = a$ and $\phi_{\theta}(x\theta(x)) = ax\theta(x)$. Bijectivity of $\phi_{\theta}$ implies $x\theta(x) \neq a^{2}$. Thus, if $x\theta(x) \neq a^{2}$, then $\theta$ is an orthomorphism given by $(ii)$. Clearly, we have 16 orthomorphisms of this form in $Orth(G)$.
	 
	 \textbf{Case(ii):} Assume $\theta(ax\theta(x)) = x$. 
	 Clearly, $\theta(a\theta(x)) = x\theta(x)$. Then $\phi_{\theta}(ax\theta(x)) = a^{-1}\theta(x)$ and $\phi_{\theta}$ $(a\theta(x)) = a^{-1}x$.
	 
	 \textbf{Subcase(a):} Assume $\theta(\theta(x)) = a$. Then $\theta(x\theta(x)) = a\theta(x)$, $\phi_{\theta}(\theta(x)) = a\theta(x)$ and $\phi_{\theta}(x\theta(x)) = ax$. Bijectivity of $\phi_{\theta}$ implies, $x\theta(x) \neq a^{2}$. Thus, if  $x\theta(x) \neq a^{2}$, then $\theta$ is an orthomorphism given by $(iii)$. Clearly, we have 16 orthomorphisms of this form in $Orth(G)$.
	 
	 \textbf{Subcase(b):} Assume $\theta(\theta(x)) = a\theta(x)$. Then $\theta(x\theta(x)) = a$, $\phi_{\theta}(\theta(x))$ = $a$ and $\phi_{\theta}(x\theta(x)) = ax\theta(x)$. Bijectivity of $\phi_{\theta}$ implies $x\theta(x) = a^{2}$. Thus, if $x\theta(x) = a^{2}$, then $\theta$ is an orthomorphism given by $(iv)$. Clearly, we have 8 orthomorphisms of this form in $Orth(G)$.
	 Hence, $\vert Orth(G) \vert = 48$.
	\end{proof}
	 \section{Clique in $Orth(\mathbb{Z}_2 \times \mathbb{Z}_4)$}
	 Let $G = \mathbb{Z}_2 \times \mathbb{Z}_4$ and $A_{ij}$, $A_{ij}^{'}$ denotes the partition of $\mathbb{Z}_2 \times \mathbb{Z}_4$ with respect to $\theta_{1}$, $\theta_{2}$ respectively as defined in Lemma \ref{le1}.
	 \begin{lemma}\label{le2}
	 	If $\theta_{1}$ and $\theta_{2} \in Orth(G)$, then $\theta_{1} \perp \theta_{2}$ if and only if 
	 	\begin{enumerate}
	 		\item [(a)] $\{\theta_{1}(x)^{-1}\theta_{2}(x) \mid x \in A_{44}\cap A_{44}^{'}\} \sqcup \\ \{\theta_{1}(x)^{-1}\theta_{2}(x)\mid x \in A_{24}\cap A_{24}^{'}\} \sqcup\\ \{\theta_{1}(x)^{-1}\theta_{2}(x) \mid x \in A_{42}\cap A_{42}^{'}\} \sqcup \\ \{\theta_{1}(x)^{-1}\theta_{2}(x)\mid x \in A_{22}\cap A_{22}^{'}\}$ = $\{x \in G \mid o(x) = 2\}$.
	 		\item [(b)] $\{\theta_{1}(x)^{-1}\theta_{2}(x) \mid x \in A_{44}\cap A_{42}^{'}\} \sqcup \\ \{\theta_{1}(x)^{-1}\theta_{2}(x)\mid x \in A_{42}\cap A_{44}^{'}\} \sqcup \\ \{\theta_{1}(x)^{-1}\theta_{2}(x) \mid x \in A_{24}\cap A_{22}^{'}\} \sqcup \\ \{\theta_{1}(x)^{-1}\theta_{2}(x)\mid x \in A_{22}\cap A_{24}^{'}\}$ = $\{x \in G \mid o(x) = 4\}$.
 		\end{enumerate}
	 \end{lemma}
 \begin{proof}
 Follows from the bijectivity of map $x \mapsto \theta_{1}(x)^{-1}\theta_{2}(x)$.
 \end{proof}
\begin{proposition}\label{pr2}
	Let $\theta_{1}$,$\theta_{2} \in Orth(G)$ such that $\vert A_{44} \cap A_{44}^{'} \vert = 2$. Then $\theta_{1} \not\perp \theta_{2}$.
\end{proposition}
\begin{proof}
	Since $\vert A_{44} \cap A_{44}^{'} \vert = 2$, $\vert A_{42} \cap A_{42}^{'} \vert = 2$.\\
	So, $\vert\{\theta_{1}(x)^{-1}\theta_{2}(x) \mid x \in A_{44}\cap A_{44}^{'}\}\vert +  \vert\{\theta_{1}(x)^{-1}\theta_{2}(x)\mid x \in A_{24}\cap A_{24}^{'}\}\vert = 4 $ $> \vert \{x \in G \mid o(x) = 2\}\vert$.
	Therefore, by  Lemma \ref{le2}, $\theta_{1}\not\perp\theta_{2}$.
\end{proof}
\begin{proposition}\label{pr3}
	Let $\theta_{1},\theta_{2} \in Orth(G)$ such that $\vert A_{44} \cap A_{44}^{'} \vert = 1$. Then $\theta_{1}\not\perp \theta_{2}$. 
\end{proposition}
\begin{proof}
	Suppose $A_{44} \cap A_{44}^{'} = \{a\}$. Then by Proposition \ref{pr1}, $A_{44} =\{a, ax\}$ and $A_{44}^{'} = \{a,ax'\}$, where $A_{22} =\{x\}$ and $A_{22}^{'} = \{x'\}$. Clearly, $x \neq x'$.\\
	\textbf{Case(1):} Assume $a \in A_{44} \cap \theta_{1}(A_{24})$ and $a \in A_{44}^{'}\cap \theta_{2}(A_{24}^{'})$. Then
	$\theta_{1}(a)^{-1}\theta_{2}(a) = (ax)^{-1}ax' = xx'$.\\
	\textbf{Subcase(a):} Assume $\theta_{1}(x) = \theta_{2}(x')$. Then $A_{22} \cup A_{24} = \{x,x',$ $\theta_{1}(x)\}$. Also,
	 	$\theta_{1}(a)^{-1}\theta_{2}(a) = xx' = \theta_{1}(x)$ and
	 	$a\theta_{1}(x) \in A_{42} \cap A_{42}^{'}. $ 
	Clearly, $\theta_{1}(A_{42}) = \{x,x'\} = \theta_{2}(A_{42}^{'})$. So, $\theta_{1}(a\theta_{1}(x))^{-1}\theta_{2}(a\theta_{2}(x'))$ $= e$ or $\theta_{1}(x)$. This is a contradiction to Lemma \ref{le2}.
	Hence $\theta_{1}\not\perp \theta_{2}$.\\
		\textbf{Subcase(b):} Assume $\theta_{1}(x) \neq \theta_{2}(x')$. Since $\theta_{2}(x') \in \{x,x',\theta_{1}(x)\}$ and $\theta_{2}(x') \neq x'$, $\theta_{2}(x') = x$. So, $a\theta_{1}(x) = ax'\theta_{2}(x') \in A_{42} \cap A_{42}^{'}$ and $\theta_{1}(x) = x'\theta_{2}(x') \in A_{24} \cap A_{24}^{'}$.
	
	\begin{enumerate}[leftmargin=*]
		\item[(i)] If $\theta_{1}(a\theta_{1}(x)) = x$ and $\theta_{2}(ax'\theta_2(x')) = x'$, then $\theta_{1}(a\theta_{1}(x))^{-1}$ $\theta_{2}(ax'\theta_{2}(x')) = xx' = \theta_1(x)$. By Lemma \ref{le2}, $\theta_{1}\not\perp \theta_{2}$ .
		\item [(ii)] If $\theta_{1}(a\theta_{1}(x)) = x\theta_{1}(x)$ and $\theta_{2}(ax'\theta_2(x')) = x'\theta_{2}(x')$, then $\theta_{1}(a\theta_{1}(x))^{-1}\theta_{2}(ax'\theta_{2}(x')) = x$.
		Also $\theta_{1}(x) = x'\theta_{2}(x') \in A_{24}\cap A_{24}^{'}$.\\
		If $x\theta_{1}(x) = a^{2}$, then $x'\theta_{2}(x') \neq a^{2} $. By Theorem \ref{th1} $(iv)$ and $(ii)$, $\theta_{1}(\theta_{1}(x))^{-1}\theta_{2}(x'\theta_{2}(x')) = (a\theta_{1}(x))^{-1}a = \theta_{1}(x)$. This is a contradiction to Lemma \ref{le2}.\\
		If $x\theta_{1}(x) \neq a^{2}$ and $x'\theta_{2}(x') = a^{2}$, then by Theorem \ref{th1} $(iii)$ and $(i)$, $\theta_{1}(\theta_{1}(x))^{-1}\theta_{2}(x'\theta_{2}(x')) = a^{-1}a\theta_{2}(x')= \theta_{2}(x')= x$. This is a contradiction to Lemma \ref{le2}.\\
		If $x\theta_{1}(x) \neq a^{2}$ and $x'\theta_{2}(x') \neq a^{2}$, then by Theorem \ref{th1} $(iii)$ and $(ii)$,
	    $\theta_{1}(\theta_{1}(x))^{-1}\theta_{2}(x'\theta_{2}(x')) = a^{-1}a = e$. This is a contradiction to Lemma \ref{le2}.
	    \item[(iii)] If $\theta_{1}(a\theta_{1}(x)) = x$ and  $\theta_{2}(ax'\theta_2(x')) = x'\theta_{2}(x')$, then $\theta_{1}(a\theta_{1}(x))^{-1}$ $\theta_{2}(ax'\theta_{2}(x')) = x'$.\\
	    If $x\theta_{1}(x) = a^{2}$, then $x'\theta_{2}(x') \neq a^{2}.$ Then by Theorem \ref{th1} $(i)$ and $(ii)$, $\theta_{1}(\theta_{1}(x))^{-1}\theta_{1}(x'\theta_{2}(x')) = a^{-1}a = e$. This is a contradiction to Lemma \ref{le2}.\\
	    If $x\theta_{1}(x) \neq a^{2}$ and $x'\theta_{2}(x') = a^{2}$, then by Theorem \ref{th1} $(ii)$ and $(i)$, $\theta_{1}(\theta_{1}(x))^{-1}\theta_{2}(x'\theta_{2}(x')) = (a\theta_{1}(x))^{-1}a\theta_{2}(x')= x'$. This is a contradiction to Lemma \ref{le2}.\\
	     If $x\theta_{1}(x) \neq a^{2}$ and $x'\theta_{2}(x') \neq  a^{2}$, then by Theorem \ref{th1} $(ii)$ and $(ii)$, $\theta_{1}(\theta_{1}(x))^{-1}$ $\theta_{2}(x'\theta_{2}(x')) = (a\theta_{1}(x))^{-1}a = \theta_{1}(x)$. This is a contradiction to Lemma \ref{le2}.
	     \item [(iv)] If $\theta_{1}(a\theta_{1}(x)) = x\theta_{1}(x)$ and $\theta_{2}(ax'\theta_2(x')) = x'$, then $\theta_{1}(a\theta_{1}(x))^{-1}$ $\theta_{2}(ax'\theta_{2}(x')) = x\theta_{1}(x)x' = e$
	    This is a contradiction to Lemma \ref{le2}.
	    Hence, $\theta_{1}\not\perp \theta_{2}$ when $a \in A_{44} \cap \theta_{1}(A_{24})$ and $a \in A_{44}^{'}\cap \theta_{2}(A_{24}^{'})$.
	\end{enumerate}
 \textbf{Case(2):} If $a \in A_{44} \cap \theta_{1}(A_{24})$ and $a \in A_{44}^{'}\cap \theta_{2}(A_{44}^{'})$, then $\theta_{1}(a) = ax, \theta_{1}(ax)=ax\theta_{1}(x)$ and $\theta_{2}(ax') = a, \theta_{2}(a)=a\theta_{2}(x')$.
 Clearly, $a \in A_{44} \cap A_{44}^{'}$ and
 	$\theta_{1}(a)^{-1}\theta_{2}(a) = (ax)^{-1}a\theta_{2}(x')= x\theta_{2}(x')$.\\
 	\textbf{Subcase(a):} Assume $\theta_{1}(x) = \theta_{2}(x')$. Then $A_{22} \cup A_{24} = \{x,x',$ $\theta_{1}(x)\}$ and $\theta_{1}(a)^{-1}\theta_{2}(a) = x'$. Clearly, $a\theta_{1}(x) = a\theta_{2}(x') \in A_{42} \cap A_{42}^{'}$ and $\theta_{1}(x) = \theta_{2}(x') \in A_{24}\cap A_{24}^{'}.$
 	\begin{enumerate}[leftmargin=*]
 		\item [(i)] If $\theta_{1}(a\theta_{1}(x)) = x$ and $\theta_{2}(a\theta_{2}(x')) = x'$, then $\theta_{1}(a\theta_{1}(x))^{-1}\theta_{2}(a$ $\theta_{2}(x')) = xx' = \theta_{1}(x)$.\\
 		If $x\theta_{1}(x) = a^{2}$, then $x'\theta_{2}(x') \neq a^{2}$. By Theorem \ref{th1} $(i)$ and $(ii)$, $\theta_{1}(\theta_{1}(x))^{-1}\theta_{2}(\theta_{2}(x')) = a^{-1}a\theta_{1}(x) = \theta_{1}(x)$. This is a contradiction to Lemma \ref{le2}.\\
 		If $x\theta_{1}(x) \neq a^{2}$ and $x'\theta_{2}(x') = a^{2}$, then by Theorem \ref{th1} $(ii)$ and $(i)$, $\theta_{1}(\theta_{1}(x))^{-1}\theta_{2}(\theta_{2}(x')) = (a\theta_{1}(x))^{-1}a = \theta_{1}(x)$. This is a contradiction to Lemma \ref{le2}.\\
 		If $x\theta_{1}(x) \neq a^{2}$ and $x'\theta_{2}(x') \neq a^{2}$, then by Theorem \ref{th1} $(ii)$ and $(ii)$, $\theta_{1}(\theta_{1}(x))^{-1}\theta_{2}(\theta_{2}(x')) = (a\theta_{1}(x))^{-1}a\theta_{2}(x') = e$. This is a contradiction to Lemma \ref{le2}.
 		\item[(ii)] If $\theta_{1}(a\theta_{1}(x)) = x\theta_{1}(x)$ and $\theta_{2}(a\theta_{2}(x')) = x'\theta_{2}(x')$, then $\theta_{1}(a\theta_{1}(x))^{-1}$ $\theta_{2}(a\theta_{2}(x)) = x\theta_{1}(x)x'\theta_{2}(x') = \theta_{1}(x)$.\\
 		If $x\theta_{1}(x) = a^{2}$, then $x'\theta_{2}(x') \neq a^{2}$. Then by Theorem \ref{th1} $(iv)$ and $(iii)$, $\theta_{1}(\theta_{1}(x))^{-1}\theta_{2}(\theta_{2}(x')) = (a\theta_{1}(x))^{-1}a = \theta_{1}(x)$. This is a contradiction to Lemma \ref{le2}.\\
 		If $x\theta_{1}(x) \neq a^{2}$ and  $x'\theta_{2}(x') = a^{2}$, then by Theorem \ref{th1} $(iii)$ and $(iv)$, $\theta_{1}(\theta_{1}(x))^{-1}\theta_{2}(\theta_{2}(x')) = a^{-1}a\theta_{2}(x') = \theta_{1}(x)$. This is a contradiction to Lemma \ref{le2}.\\
 		If $x\theta_{1}(x) \neq a^{2}$ and  $x'\theta_{2}(x') \neq a^{2}$, then by Theorem \ref{th1} $(iii)$ and $(iii)$, $\theta_{1}(\theta_{1}(x))^{-1}\theta_{2}(\theta_{2}(x')) = a^{-1}a = e$. This is a contradiction to Lemma \ref{le2}.
 		\item[(iii)] If $\theta_{1}(a\theta_{1}(x)) = x$ and $\theta_{2}(a\theta_{2}(x')) = x'\theta_{2}(x')$, then $\theta_{1}(a\theta_{1}(x))^{-1}$ $\theta_{2}(a\theta_{2}(x')) = xx'\theta_{2}(x') = e$. This is a contradiction to Lemma \ref{le2}.
 		\item[(iv)] If $\theta_{1}(a\theta_{1}(x)) = x\theta_{1}(x)$ and $\theta_{2}(a\theta_{2}(x')) = x'$, then $\theta_{1}(a\theta_{1}(x))^{-1}$ $\theta_{2}(a\theta_{2}(x')) = x\theta_{1}(x)x' = e$. This is a contradiction to Lemma \ref{le2}. Hence, $\theta_{1} \not\perp \theta_{2}$.
 	 	 \end{enumerate}
  	 \textbf{Subcase(b):} Assume $\theta_{1}(x) \neq \theta_{2}(x')$. Since $\theta_{2}(x^{'}) \in \{x,x^{'},\theta_{1}(x)\}$ and $\theta_{2}(x') \neq x'$, $\theta_{2}(x') = x$. As $a \in A_{44}\cap A_{44}^{'}$, so $\theta_{1}(a)^{-1}\theta_{2}(a) = x\theta_{2}(x') = e$. This is a contradiction to Lemma \ref{le2}.
  	 Thus, $\theta_{1} \not\perp \theta_{2}$ when $a \in A_{44} \cap \theta_{1}(A_{24})$ and $a \in A_{44}^{'}\cap \theta_{2}(A_{44}^{'})$.
\end{proof}
\begin{proposition}\label{pr4}
	Let $\theta_{1}, \theta_{2} \in Orth(G)$ and $\vert A_{44}\cap A_{44}^{'}\vert = 0$. Then $\theta_{1} \perp \theta_{2}.$
\end{proposition}
\begin{proof}
	If $\vert A_{44}\cap A_{44}^{'} \vert = 0$ then $\vert A_{22}\cap A_{22}^{'} \vert = 1$ and $\vert A_{24}\cap A_{24}^{'} \vert = 2$. Also $\vert A_{44}\cap A_{42}^{'} \vert = 2$ and $\vert A_{42}\cap A_{44}^{'} \vert = 2$.
	Consider the orthomorphism of the form $\theta_{1} = (a, ax, ax\theta_{1}(x), x\theta_{1}(x),a\theta_{1}(x),x,\theta_{1}(x))$ where $x\theta_{1}(x) = a^{2}$.
	Here, $A_{22} = \{x\}$, $A_{44} = \{a,  ax\}$ , $A_{42} = \{ax\theta_{1}(x), a\theta_{1}(x)\}$ and $A_{24} = \{\theta_{1}(x), x\theta_{1}(x)\}$. If $\theta_{2}$ is orthogonal to $\theta_{1}$ then, $A_{22}^{'} = \{x\}$, $A_{44}^{'} =\{ax\theta_{1}(x),$ $ a\theta_{1}(x)\}$, $A_{42}^{'} =  \{a,  ax\}$ and $A_{24}^{'} = \{\theta_{1}(x), x\theta_{1}(x)\}$.\\
	\textbf{Case(1):} Assume  $a\theta_{1}(x) \in A_{44}^{'} \cap \theta_{2}(A_{24}^{'})$.
	 Then by Proposition \ref{pr1}, $\theta_{2}(a\theta_{1}(x)) = ax\theta_{1}(x)$ and $\theta_{2}(ax\theta_{1}(x)) = a$ as $\theta_{2}(x) = x\theta_{1}(x).$ Since $\theta_{1}(x\theta_{1}(x))$ $= a\theta_{1}(x)$, $\theta_{2}(x\theta_{1}(x))= ax$ and $\theta_{2}(\theta_{1}(x)) = a\theta_{1}(x)$. Then $\phi_{\theta_{2}}(a\theta_{1}(x)) = x ,\phi_{\theta_{2}}(ax\theta_{1}(x)) = x\theta_{1}(x), \phi_{\theta_{2}}(x) = \theta_{1}(x), \phi_{\theta_{2}}(x\theta_{1}
(x))= a\theta_{1}(x), \phi_{\theta_{2}}(\theta_{1}(x)) = a$.\\
\textbf{Subcase(a):} Assume $\theta_{2}(a) = x$. Then $\theta_{2}(ax) = \theta_{1}(x)$, $\phi_{\theta_{2}}(a) = a^{-1}x$ and   $\phi_{\theta_{2}}(ax) = (ax)^{-1}\theta_{1}(x) = a$, which is not an orthomorphism.\\
\textbf{Subcase(b):} Assume $\theta_{2}(a) = \theta_{1}(x)$. Then $\theta_{2}(ax) = x$, $\phi_{\theta_{2}}(a) = a^{-1}\theta_{1}(x) = ax$ and   $\phi_{\theta_{2}}(ax) = (ax)^{-1}x = a^{-1}$.
In this case $\theta_{2}$ becomes an orthomorphism given by
$\theta_{2} = (a\theta_{1}(x),ax\theta_{1}(x),a, \theta_{1}(x))$ $(x,x\theta_{1}(x),ax)$.\\
Now,
$$\theta_{1}(y)^{-1}\theta_{2}(y)=\begin{dcases}
	\theta_{1}(x)x\theta_{1}(x) = x & y = x \in A_{22} \cap A_{22}^{'}\\
	a^{-1}a\theta_{1}(x) = \theta_{1}(x) & y = \theta_{1}(x) \in A_{24} \cap A_{24}^{'}\\
	a^{-1}\theta_{1}(x)ax = x\theta_{1}(x) & y = x\theta_{1}(x) \in A_{24} \cap A_{24}^{'}\\
	a^{-1}x\theta_{1}(x) = a & y = a \in A_{44} \cap A_{42}^{'}\\
	a^{-1}x\theta_{1}(x)x = ax & y = ax \in A_{44} \cap A_{42}^{'}\\
	xax\theta_{1}(x) = a\theta_{1}(x) & y = a\theta_{1}(x) \in A_{42} \cap A_{44}^{'}\\
	x\theta_{1}(x)a = ax\theta_{1}(x) & y = ax\theta_{1}(x) \in A_{42} \cap A_{44}^{'}\\
\end{dcases}$$
Clearly, $y \mapsto \theta_{1}(y)^{-1}\theta_{2}(y)$ is a bijective map. Thus, $\theta_{1} \perp \theta_{2}$.\\
\textbf{Case(2):} Assume $a\theta_{1}(x) \in A_{44}^{'} \cap \theta_{2}(A_{44}^{'}).$
 Then by Proposition \ref{pr1}, $\theta_{2}(a\theta_{1}(x)) = ax, \theta_{2}(ax\theta_{1}(x)) = a\theta_{1}(x)$ as $\theta_{2}(x) = x\theta_{1}(x).$\\
Since $\theta_{1}(\theta_{1}(x)) = a$, $\theta_{2}(\theta_{1}(x)) = ax\theta_{1}(x)$ and $\theta_{2}(x\theta_{1}(x)) = a$. Then
$\phi_{\theta_{2}}(a\theta_{1}(x)) = x\theta_{1}(x),\phi_{\theta_{2}}(ax\theta_{1}(x)) = x, \phi_{\theta_{2}}(x) = \theta_{1}(x),$ $\phi_{\theta_{2}}(x\theta_{1}
(x))= ax\theta_{1}(x) = a^{-1}, \phi_{\theta_{2}}(\theta_{1}(x)) = ax$.\\
\textbf{Subcase(a):}  Assume $\theta_{2}(a) = \theta_{1}(x)$. Then $\theta_{2}(ax) = x$,
$\phi_{\theta_{2}}(a) = a^{-1}\theta_{1}(x)$ and   $\phi_{\theta_{2}}(ax) = (ax)^{-1}x = a^{-1}$ which is not an orthomorphism.\\
\textbf{Subcase(b):} Assume $\theta_{2}(a) = x$. Then $\theta_{2}(ax) = \theta_{1}(x)$, 
$\phi_{\theta_{2}}(a) = a^{-1}x$ and   $\phi_{\theta_{2}}(ax) = (ax)^{-1}\theta_{1}(x)$.
In this case $\theta_{2}$ becomes an orthomorphism given by $\theta_{2}=(ax\theta_{1}(x),a\theta_{1}(x),ax,\theta_{1}(x))$ $(x,x\theta_{1}(x),a)$.\\
Now,
$$\theta_{1}(y)^{-1}\theta_{2}(y)=\begin{dcases}
	\theta_{1}(x)x\theta_{1}(x) = x & y = x \in A_{22} \cap A_{22}^{'}\\
	a^{-1}ax\theta_{1}(x) = x\theta_{1}(x) & y = \theta_{1}(x) \in A_{24} \cap A_{24}^{'}\\
	a^{-1}\theta_{1}(x)a = \theta_{1}(x) & y = x\theta_{1}(x) \in A_{24} \cap A_{24}^{'}\\
	a^{-1}xx = ax\theta_{1}(x) & y = a \in A_{44} \cap A_{42}^{'}\\
	a^{-1}x\theta_{1}(x)\theta_{1}(x) = a\theta_{1}(x) & y = ax \in A_{44} \cap A_{42}^{'}\\
	xax = a & y = a\theta_{1}(x) \in A_{42} \cap A_{44}^{'}\\
	x\theta_{1}(x)a\theta_{1}(x) = ax & y = ax\theta_{1}(x) \in A_{42} \cap A_{44}^{'}\\
\end{dcases}$$
Clearly, $y \mapsto \theta_{1}(y)^{-1}\theta_{2}(y)$ is a bijective map. Thus, $\theta_{1} \perp \theta_{2}$.\\
Hence, If $\theta_{1}=(a,ax, ax\theta_{1}(x),x\theta_{1}(x),a\theta_{1}(x),x,\theta_{1}(x))$ where $x\theta_{1}(x) = a^{2}$ then $\theta_{1}$ is orthogonal to 
$\theta_{2}=(a\theta_{1}(x),ax\theta_{1}(x),a,\theta_{1}(x))(x,x\theta_{1}(x),ax)$
	and $\theta_{3}=(ax\theta_{1}(x),a\theta_{1}(x),ax,\theta_{1}(x))(x,x\theta_{1}(x),a)$.
\end{proof}
	Similarly, calculating the other cases, the following Table 1 has been constructed:
	\begin{center}
			\begin{table}[h!]
				\renewcommand{\arraystretch}{1.2}
		\begin{tabular}{|c|c|}
			\hline
			$\theta_{1}$ & $\theta_{2}$, $\theta_{3}$ \\ \hline
			\shortstack{$(a,ax, ax\theta_{1}(x),x\theta_{1}(x),a\theta_{1}(x),x,\theta_{1}(x))$\\  where $x\theta_{1}(x) = a^{2}$} & \shortstack{$(a\theta_{1}(x),ax\theta_{1}(x),a,\theta_{1}(x))(x,x\theta_{1}(x),ax)$,\\ $(ax\theta_{1}(x),a\theta_{1}(x),ax,\theta_{1}(x))(x,x\theta_{1}(x),a)$}.\\ \hline
			\shortstack{$(a,ax, ax\theta_{1}(x),x,\theta_{1}(x),a\theta_{1}(x),x\theta_{1}(x))$ \\ where $x\theta_{1}(x) = a^{2}$} &  \shortstack{$(a\theta_{1}(x),ax\theta_{1}(x),a,x,x\theta_{1}(x))(ax,\theta_{1}(x))$, \\ $(ax\theta_{1}(x),a\theta_{1}(x),ax,x,x\theta_{1}(x))(a,\theta_{1}(x))$}.\\ \hline
			\shortstack{$(a, ax, ax\theta_{1}(x), x\theta_{1}(x))(\theta_{1}(x), a\theta_{1}(x),x)$\\where $x\theta_{1}(x) \neq a^{2}$ and  $x = a^{2}$}	& \shortstack{$(a\theta_{1}(x),ax\theta_{1}(x),a,x,x\theta_{1}(x))(ax,\theta_{1}(x))$,\\ $(ax\theta_{1}(x),a\theta_{1}(x),ax,x,x\theta_{1}(x))(a,\theta_{1}(x))$.}\\ \hline
			\shortstack{$(a, ax, ax\theta_{1}(x), x\theta_{1}(x))(\theta_{1}(x), a\theta_{1}(x),x)$\\where $x\theta_{1}(x) \neq a^{2}$ and  $\theta_{1}(x) = a^{2}$} & \shortstack{$(a\theta_{1}(x),ax\theta_{1}(x),a,\theta_{1}(x),ax,x,x\theta_{1}(x))$, \\$(ax\theta_{1}(x),a\theta_{1}(x),ax,\theta_{1}(x),a,x,x\theta_{1}(x))$}. \\ \hline
			\shortstack{$(a, ax, ax\theta_{1}(x),x,\theta_{1}(x)) (a\theta_{1}(x),x\theta_{1}(x))$ \\ where $x\theta_{1}(x) \neq a^{2}$ and $x = a^{2}$} & \shortstack{$(a\theta_{1}(x),ax\theta_{1}(x),a,\theta_{1}(x))(ax,x,x\theta_{1}(x))$, \\$(ax\theta_{1}(x),a\theta_{1}(x),ax,\theta_{1}(x))(a,x,x\theta_{1}(x))$}. \\ \hline
			\shortstack{$(a, ax, ax\theta_{1}(x),x,\theta_{1}(x)) (a\theta_{1}(x),x\theta_{1}(x))$ \\ where $x\theta_{1}(x) \neq a^{2}$ and $\theta_{1}(x) = a^{2}$} & \shortstack{$(a\theta_{1}(x),ax\theta_{1}(x),a,x,x\theta_{1}(x),ax,\theta_{1}(x))$,
				\\$(ax\theta_{1}(x),a\theta_{1}(x),ax,x,x\theta_{1}(x),a,\theta_{1}(x))$}.\\ \hline
		\end{tabular}
	\caption{ $\theta_{1} \perp \theta_{2}$ and $\theta_{1} \perp \theta_{3}$}
	\end{table}
	\end{center}
\begin{corollary}\label{co3}
	$\omega(\mathbb{Z}_{2}\times \mathbb{Z}_{4}) = 2$.
\end{corollary}
\begin{proof}
	Clearly, there are two orthomorphism orthogonal to a given orthomorphism and they cannot be orthogonal to each other as their $A_{44}$ are same. Hence, $\omega(\mathbb{Z}_{2}\times \mathbb{Z}_{4}) = 2$.
\end{proof}
\begin{corollary}
	\begin{enumerate}
		\item [(i)] Two orthomorphism $\psi_{1}$ and $\psi_{2}$ which are orthogonal to $\theta$ are conjugate to each other by an element $\alpha = (a,ax)(ax\theta(x), a\theta(x))$ in $Aut(Z_{2} \times Z_{4})$ where $A_{44} = \{a, ax\}$ and $A_{22} = \{x\}$ of $\theta$. Also $\psi_{1}$ and $\psi_{2}$ are also orthogonal to $\alpha\theta\alpha^{-1} = \theta^{\alpha}$.
		\item[(ii)] $Orth(Z_{2} \times Z_{4})$ consists of 12 disjoint 4-cycles. Each 4-cycle is given by Figure 2.
		\begin{figure}[H]
			\begin{center}
				\includegraphics[scale = 1]{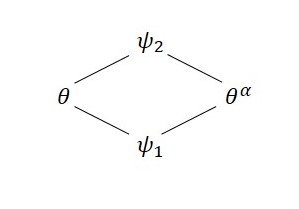} \caption{}
			\end{center}
			
		\end{figure}
	\end{enumerate}
\end{corollary}

\end{document}